\documentclass[12pt]{article}%
\usepackage{amsmath}
\usepackage{amsfonts}
\usepackage{amssymb}
\usepackage{graphicx}%
\newtheorem{theorem}{Theorem}

\newtheorem{lemma}[theorem]{Lemma}

\newenvironment{proof}[1][Proof]{\noindent\textbf{#1.} }{\ \rule{0.5em}{0.5em}}
\begin{document}

\title{ Variations on a Theme by James Stirling}
\author{Diego Dominici \thanks{e-mail: dominicd@newpaltz.edu}\\Department of Mathematics\\State University of New York at New Paltz\\75 S. Manheim Blvd. Suite 9\\New Paltz, NY 12561-2443\\USA\\Phone: (845) 257-2607\\Fax: (845) 257-3571 }
\date{\textit{Dedicated to Donald Silberger on the occasion of his 76th birthday} }
\maketitle

\section{Introduction}

The most widely known and used result in asymptotics is probably
\textit{Stirling's formula},
\begin{equation}
n!\sim\sqrt{2\pi n}n^{n}e^{-n},\quad n\rightarrow\infty\label{n!}%
\end{equation}
named after James Stirling (May 1692 -- 5 Dec. 1770). The formula provides an
extremely accurate approximation of the factorial numbers for large values of
$n.$ The asymptotic formula (\ref{n!}), for which Stirling is best known,
appeared as Example 2 to Proposition 28 in his most important work,
\textit{Methodus Differentialis, }published in 1730 \cite{MR1996414}. In it,
he asserted that $\log\left(  n!\right)  $ is approximated by
\textquotedblleft three or four terms" of the series%
\begin{align}
&  \left(  n+\frac{1}{2}\right)  \log\left(  n+\frac{1}{2}\right)  -a\left(
n+\frac{1}{2}\right)  +\frac{1}{2}\log\left(  2\pi\right) \label{stirl1}\\
&  -\frac{a}{24\left(  n+\frac{1}{2}\right)  }+\frac{7a}{2880\left(
n+\frac{1}{2}\right)  ^{3}}-\cdots,\nonumber
\end{align}
where $\log$ means the base-10 logarithm and $a=\left[  \ln(10)\right]  ^{-1}$
\cite{MR747189}.

In 1730 Stirling wrote to Abraham De Moivre (26 May 1667 -- 27 Nov. 1754)
pointing out some errors that he had made in a table of logarithms of
factorials in the book and also telling him about (\ref{stirl1}). After seeing
Stirling's results, De Moivre derived the formula%
\begin{equation}
\ln\left[  \left(  n-1\right)  !\right]  \sim\left(  n-\frac{1}{2}\right)
\ln(n)-n+\frac{1}{2}\ln\left(  2\pi\right)  +%
{\displaystyle\sum\limits_{k\geq1}}
\frac{B_{2k}}{2k\left(  2k-1\right)  n^{2k-1}}, \label{n!1}%
\end{equation}
which he published in his \textit{Miscellaneis Analyticis Supplementum} a few
months later. Equation (\ref{n!1}) is called \textit{Stirling's series }and
the numbers $B_{k}$ are called the \textit{Bernoulli numbers}, and are defined
by
\begin{equation}
B_{0}=1,\quad%
{\displaystyle\sum\limits_{j=0}^{k}}
\binom{k+1}{j}B_{j}=0,\quad k\geq1. \label{bernoulli}%
\end{equation}
Clearly Stirling and De Moivre regularly corresponded around this time, for in
September 1730 Stirling relates the new results of De Moivre in a letter to
Gabriel Cramer.

In 1729 Leonhard Euler (15 April 1707 -- 18 Sept. 1783) proposed a
generalization of the factorial function from natural numbers to positive real
numbers \cite{MR0106810}. It is called the \textit{gamma function,}
$\Gamma(z)$, which he defined as
\begin{equation}
\Gamma(z)=\underset{n\rightarrow\infty}{\lim}\frac{n!n^{z}}{z(z+1)\cdots
(z+n)}, \label{euler}%
\end{equation}
and it is related to the factorial numbers by%
\[
\Gamma\left(  n+1\right)  =n!,\quad n=0,1,2,\ldots.
\]
From Euler's definition (\ref{euler}), we immediately obtain the fundamental
relation%
\begin{equation}
\Gamma(z+1)=z\Gamma(z) \label{difference}%
\end{equation}
and the value $\Gamma(1)=1.$ In fact, the gamma function is completely
characterized by the \textit{Bohr-Mollerup theorem} \cite{Borh}:

\begin{theorem}
The gamma function is the only function $\Gamma:\left(  0,\infty\right)
\rightarrow\left(  0,\infty\right)  $ which satisfies

\begin{enumerate}
\item $\Gamma(1)=1$

\item $\Gamma(x+1)=x\Gamma(x)$

\item $\ln\left[  \Gamma(x)\right]  $ is convex

for all $x\in\left(  0,\infty\right)  .$
\end{enumerate}
\end{theorem}

\begin{proof}
See \cite{MR0165148} and \cite{MR0492418}.
\end{proof}

Another complex-analytic characterization is due to Wielandt \cite{MR1430098}:

\begin{theorem}
The gamma function is the only holomorphic function in the right half plane
$\mathbb{A}$ satisfying

\begin{enumerate}
\item $\Gamma(1)=1$

\item $\Gamma(z+1)=z\Gamma(z)$ for all $z\in\mathbb{A}$

\item $\Gamma\left(  z\right)  $ is bounded in the strip $1\leq
\operatorname{Re}(z)<2$
\end{enumerate}
\end{theorem}

\begin{proof}
See \cite{MR1376175}.
\end{proof}

In terms of $\Gamma(z),$ we can re-write (\ref{n!}) as
\begin{equation}
\ln\left[  \Gamma(z)\right]  \sim\mathcal{P}(z),\quad z\rightarrow
\infty\label{Gamma1}%
\end{equation}
with%
\begin{equation}
\mathcal{P}(z)=z\ln\left(  z\right)  -z-\frac{1}{2}\ln\left(  z\right)
+\frac{1}{2}\ln\left(  2\pi\right)  \label{P}%
\end{equation}
and (\ref{n!1}) in the form%
\begin{equation}
\ln\left[  \Gamma(z)\right]  \sim\mathcal{P}(z)+R_{N}(z),\quad z\rightarrow
\infty\label{Gamma2}%
\end{equation}
where $R_{0}(z)=0$ and%
\begin{equation}
R_{N}(z)=%
{\displaystyle\sum\limits_{k=1}^{N}}
\frac{B_{2k}}{2k\left(  2k-1\right)  z^{2k-1}},\quad N\geq1. \label{R}%
\end{equation}
Estimations of the remainder $\ln\left[  \Gamma(z)\right]  -\mathcal{P}%
(z)-R_{N}(z)$ were computed in \cite{MR0295539}.

\section{Previous results}

Over the years, there have been \textbf{many} different approaches to the
derivation of (\ref{Gamma1}) and (\ref{Gamma2}), including:

\begin{enumerate}
\item Aissen \cite{MR0064914} studied the sequence $V_{n}=\frac{n^{n}e^{-n}%
}{n!}.$ Using his lemma

\begin{lemma}
If
\[
\frac{y_{n+1}}{y_{n}}=1+\frac{\alpha}{n}+O\left(  n^{-2}\right)
\]
and $y_{n}\neq0$ for all $n,$ then $y_{n}\sim Cn^{\alpha},\quad n\rightarrow
\infty$ for some non-zero constant $C.$
\end{lemma}

he showed that $n!\sim C\sqrt{n}n^{n}e^{-n}.$

\item Bender \& Orszag \cite{MR538168}, Bleistein \& Handelsman
\cite{MR863284}, Diaconis \& Freedman \cite{MR827588}, Dingle \cite{MR0499926}%
, Olver \cite{MR1429619} and Wong \cite{MR1851050} applied Laplace's method to
the \textit{Euler integral of the second kind}%
\begin{equation}
\Gamma\left(  z\right)  =%
{\displaystyle\int\limits_{0}^{\infty}}
t^{z-1}e^{-t}dt,\quad\operatorname{Re}(z)>0. \label{eulerint}%
\end{equation}

\item Bender \& Orszag \cite{MR538168} and Temme \cite{MR1376370} used
\textit{Hankel's contour integral} \cite{MR1225604}
\[
\frac{1}{\Gamma(z)}=\frac{1}{2\pi\mathrm{i}}%
{\displaystyle\int\limits_{-\infty}^{\left(  0+\right)  }}
t^{-z}e^{t}dt
\]
and the method of steepest descent.

\item Bleistein \& Handelsman \cite{MR863284}, Lebedev \cite{MR0350075},
Sasv\'{a}ri \cite{MR1671869} and Temme \cite{MR1376370} used \textit{Binet's
first formula}%
\[
\ln\left[  \Gamma(z)\right]  =\mathcal{P}(z)+%
{\displaystyle\int\limits_{0}^{\infty}}
\frac{1}{t}\left(  \frac{1}{2}-\frac{1}{t}+\frac{1}{e^{t}-1}\right)
e^{-tz}dt,\quad\operatorname{Re}(z)>0
\]
and
\[
\frac{1}{t}\left(  \frac{1}{2}-\frac{1}{t}+\frac{1}{e^{t}-1}\right)  =%
{\displaystyle\sum\limits_{k\geq1}}
\frac{B_{2k}}{\left(  2k\right)  !}t^{2k-2},\quad\left\vert t\right\vert
<2\pi.
\]

\item Blyth \& Pathak \cite{MR1540867} and Khan \cite{MR0339316} used
probabilistic arguments, applying the Central Limit Theorem and the limit
theorem for moment generating functions to Gamma and Poison random variables.

\item Coleman \cite{MR1527865} defined%
\[
c_{n}=\left(  n+\frac{1}{2}\right)  \ln(n)-n+1-\ln\left(  n!\right)  ,
\]
and showed that $c_{n}\rightarrow1-\frac{1}{2}\ln\left(  2\pi\right)  $ as
$n\rightarrow\infty.$ A similar result was proved by Aissen \cite{MR0064914},
using the concavity of $\ln(x).$

\item Dingle \cite{MR0499926} used \textit{Weierstrass' infinite product}%
\[
\frac{1}{\Gamma(z)}=ze^{\gamma z}%
{\displaystyle\prod\limits_{n=1}^{\infty}}
\left[  \left(  1+\frac{z}{n}\right)  e^{-\frac{z}{n}}\right]  ,
\]
(where $\gamma$ is \textit{Euler's constant)} and Mellin transforms.

\item Feller \cite{MR0222535}, \cite{MR0230020} proved the identity%
\begin{equation}
\ln\left(  n!\right)  -\frac{1}{2}\ln(n)=I(n)-I\left(  \frac{1}{2}\right)  +%
{\displaystyle\sum\limits_{k=1}^{n-1}}
\left(  a_{k}-b_{k}\right)  +a_{n}, \label{feller}%
\end{equation}
where%
\[
I(n)=%
{\displaystyle\int\limits_{0}^{n}}
\ln(t)dt,\quad a_{k}=%
{\displaystyle\int\limits_{k-\frac{1}{2}}^{k}}
\ln\left(  \frac{k}{t}\right)  dt,\quad b_{k}=%
{\displaystyle\int\limits_{k}^{k+\frac{1}{2}}}
\ln\left(  \frac{t}{k}\right)  dt
\]
and showed that%
\[%
{\displaystyle\sum\limits_{k=1}^{\infty}}
\left(  a_{k}-b_{k}\right)  -I\left(  \frac{1}{2}\right)  =\frac{1}{2}%
\ln\left(  2\pi\right)  .
\]
There is a mistake in his Equation (2.4), where he states that%
\[
\ln\left(  n!\right)  -\frac{1}{2}\ln(n)+I(n)-I\left(  \frac{1}{2}\right)  =%
{\displaystyle\sum\limits_{k=1}^{n-1}}
\left(  a_{k}-b_{k}\right)  +a_{n},
\]
instead of (\ref{feller}).

\item Hayman \cite{MR0080749} used the exponential generating function%
\[
e^{z}=%
{\displaystyle\sum\limits_{k=0}^{\infty}}
\frac{1}{k!}z^{k}%
\]
and his method for admissible functions.

\item Hummel \cite{MR1524780}, established the inequalities%
\[
\frac{11}{12}<r_{n}+\frac{1}{2}\ln\left(  2\pi\right)  <1,\quad n=2,3,\ldots,
\]
where%
\[
r_{n}=\ln\left(  \frac{n!e^{n}}{\sqrt{2\pi n}n^{n}}\right)  .
\]
Impens \cite{MR2182667}, \cite{MR2024001}, showed that for $x>0$%
\[
R_{2n}(x)<\ln\left[  \Gamma(x)\right]  -\mathcal{P}(x)<R_{2m+1}(x),\quad
n,m\geq0,
\]
where $R_{n}(x)$ was defined in (\ref{R}). Maria \cite{MR0200485} showed that%
\[
\left[  12n+\frac{3}{2\left(  2n+1\right)  }\right]  ^{-1}<r_{n},\quad
n=1,2,\ldots.
\]
Mermin \cite{MR740240} proved the identity%
\[
e^{r_{n}}=%
{\displaystyle\prod\limits_{k=n}^{\infty}}
e^{-1}\left(  1+\frac{1}{k}\right)  ^{k+\frac{1}{2}},
\]
which he used to show that $r_{n}\sim$ $R_{3}(n).$ Michel \cite{MR1903427},
proved the inequality%
\[
\left\vert e^{r_{n}}-1-\frac{1}{12n}-\frac{1}{288n^{2}}\right\vert \leq
\frac{1}{360n^{3}}+\frac{1}{108n^{4}},\quad n=3,4\ldots.
\]
Nanjundiah \cite{MR0117359}, showed that%
\[
R_{2}(n)<r_{n}<R_{1}(n)=1,2,\ldots.
\]
Robbins \cite{MR0069328}, established the double inequality%
\[
\frac{1}{12n+1}<r_{n}<\frac{1}{12n},\quad n=1,2,\ldots.
\]

\item Marsaglia \& Marsaglia \cite{MR1080390} derived from (\ref{eulerint})
the asymptotic expansion%
\[
n!\sim n^{n+1}e^{-n}%
{\displaystyle\sum\limits_{k=1}^{\infty}}
b_{k}\left(  \frac{2}{n}\right)  ^{\frac{k}{2}}\Gamma\left(  \frac{k}%
{2}\right)  k,
\]
where the generating function $G(z)=%
{\displaystyle\sum\limits_{k\geq0}}
b_{k}z^{k}$ is defined by%
\[
G(z)\exp\left[  1-G(z)\right]  =\exp\left(  -\frac{1}{2}z^{2}\right)  ,\quad
G^{\prime}(0)=1.
\]

\item Namias \cite{MR824587} introduced the function $F(n)=\frac{\Gamma
(n)}{\mathcal{P}(n)},$ with $\mathcal{P}(n)$ defined in (\ref{P}). From
\textit{Legendre's duplication formula}%
\begin{equation}
\Gamma(2n)=\frac{2^{2n-1}}{\sqrt{\pi}}\Gamma(n)\Gamma\left(  n+\frac{1}%
{2}\right)  , \label{dupli}%
\end{equation}
he derived a functional equation for $F(n)$%
\[
\frac{F(2n)}{F(n)F\left(  n-\frac{1}{2}\right)  }=\sqrt{e}\left(  1-\frac
{1}{2n}\right)  ^{n},
\]
from which he obtained (\ref{Gamma2}). He also considered the triplication
case, using \textit{Gauss' multiplication formula}%
\[
\Gamma(mz)=\left(  \sqrt{2\pi}\right)  ^{1-m}m^{mz-\frac{1}{2}}%
{\displaystyle\prod\limits_{k=0}^{m-1}}
\Gamma\left(  z+\frac{k}{m}\right)  ,\quad m=2,3,\ldots
\]
with $m=3.$ His results where extended by Deeba \& Rodriguez in
\cite{MR1104307}.

\item Olver \cite{MR1429619} used Euler's definition (\ref{euler})
\[
\ln\left[  \Gamma(z)\right]  =-\ln\left(  z\right)  +\underset{n\rightarrow
\infty}{\lim}z\ln\left(  n\right)  +%
{\displaystyle\sum\limits_{k=1}^{n}}
\left[  \ln\left(  k\right)  -\ln\left(  z+k\right)  \right]
\]
and the Euler-Maclaurin formula. A similar analysis was done by Knopp
\cite{knopp} and Wilf \cite{MR542284}.

\item Patin \cite{MR979596} used (\ref{eulerint}) and the Lebesgue Dominated
Convergence Theorem.

\item Whittaker \& Watson \cite{MR1424469} used \textit{Binet's second
formula}
\[
\ln\left[  \Gamma(z)\right]  =\mathcal{P}(z)+2%
{\displaystyle\int\limits_{0}^{\infty}}
\frac{\arctan\left(  \frac{t}{z}\right)  }{e^{2\pi t}-1}dt,\quad
\operatorname{Re}(z)>0
\]
and
\[
\arctan(x)=x%
{\displaystyle\sum\limits_{k=0}^{\infty}}
\frac{\left(  -1\right)  ^{k}}{2k+1}x^{2k},\quad\left\vert x\right\vert
\leq1.
\]

\end{enumerate}

Thus, there have been a huge variety of approaches to Stirling's result,
ranging from elementary to heavy-machinery methods. In an effort to join such
illustrious company, we present still another direction for deriving
(\ref{Gamma2}).

Our starting point shall be the difference equation (\ref{difference}). A
parallel approach was considered in \cite{askey}. For a different analysis of
(\ref{difference}) using the method of controlling factors, see
\cite{MR538168}. Extensions and other applications of the method used can be
found in \cite{MR1373150}, \cite{MR0225511} and \cite{dom}.

\section{Asymptotic analysis}

\subsection{Stirling's formula\label{section}}

We begin with a derivation of (\ref{Gamma1}), to better illustrate how the
method works. We assume that%
\begin{equation}
\ln\left[  \Gamma(z)\right]  \sim f(z)+g(z),\quad z\rightarrow\infty
\label{asymp1}%
\end{equation}
with
\begin{equation}
g=o(f),\quad z\rightarrow\infty. \label{o}%
\end{equation}
Using (\ref{asymp1}) in (\ref{difference}), we have%
\begin{equation}
f(z+1)-f(z)+g(z+1)-g(z)\sim\ln(z). \label{asymp2}%
\end{equation}
Expanding $f(z+1)$ and $g(z+1)$ in a Taylor series, we obtain%
\begin{equation}
f^{\prime}(z)+\frac{1}{2}f^{\prime\prime}(z)+g^{\prime}(z)\sim\ln(z).
\label{asymp3}%
\end{equation}
From (\ref{o}) and (\ref{asymp3}) we get the system%
\[
f^{\prime}(z)=\ln(z),\quad\frac{1}{2}f^{\prime\prime}(z)+g^{\prime}(z)=0
\]
and thus,%
\begin{equation}
f(z)=z\ln(z)-z,\quad g(z)=-\frac{1}{2}\ln(z)+C. \label{asymp4}%
\end{equation}
To find the constant $C$ in (\ref{asymp4}), we replace $\Gamma(z)\sim
e^{f(z)+g(z)}$ in (\ref{dupli}) and obtain%
\[
e^{f(2z)+g(2z)}\sim\frac{2^{2z-1}}{\sqrt{\pi}}e^{f(z)+g(z)}e^{f\left(
z+\frac{1}{2}\right)  +g\left(  z+\frac{1}{2}\right)  },
\]
or%
\[
e^{C-\frac{1}{2}}\left(  1+\frac{1}{2z}\right)  ^{z}\sim\sqrt{2\pi},\quad
z\rightarrow\infty,
\]
from which we conclude that $C=\ln\left(  \sqrt{2\pi}\right)  .$

Hence, we have shown that%
\[
\ln\left[  \Gamma(z)\right]  \sim z\ln(z)-z-\frac{1}{2}\ln(z)+\frac{1}{2}%
\ln\left(  2\pi\right)  ,\quad z\rightarrow\infty.
\]

\subsection{Stirling's series}

To extend the result of the previous section, we now assume that%

\begin{equation}
\ln\left[  \Gamma(z)\right]  \sim%
{\displaystyle\sum\limits_{k=0}^{N}}
f_{k}(z),\quad z\rightarrow\infty\label{ser1}%
\end{equation}
with%
\begin{equation}
f_{k+1}=o\left(  f_{k}\right)  ,\quad z\rightarrow\infty,\quad k=0,1,\ldots
,N-1. \label{o1}%
\end{equation}
Using (\ref{ser1}) in (\ref{difference}) we have%
\begin{equation}%
{\displaystyle\sum\limits_{k=0}^{N}}
f_{k}(z+1)-f_{k}(z)\sim\ln\left(  z\right)  ,\quad z\rightarrow\infty.
\label{ser2}%
\end{equation}
Replacing the Taylor series
\[
f_{k}(z+1)=%
{\displaystyle\sum\limits_{j\geq0}}
\frac{1}{j!}\frac{d^{j}}{dz^{j}}f_{k}(z)
\]
in (\ref{ser2}), we have
\begin{equation}%
{\displaystyle\sum\limits_{k=0}^{N}}
{\displaystyle\sum\limits_{j\geq1}}
\frac{1}{j!}\frac{d^{j}}{dz^{j}}f_{k}(z)\sim\ln\left(  z\right)  ,\quad
z\rightarrow\infty. \label{ser3}%
\end{equation}

From (\ref{o1}), we obtain the system of ODEs%
\[
\frac{d}{dz}f_{0}=\ln\left(  z\right)
\]
and%
\[%
{\displaystyle\sum\limits_{j=0}^{k-1}}
\frac{1}{\left(  k+1-j\right)  !}\frac{d^{k+1-j}}{dz^{k+1-j}}f_{j}(z)+\frac
{d}{dz}f_{k}=0,\quad k\geq1,
\]
which imply%
\begin{equation}
f_{0}(z)=z\ln(z)-z \label{f0}%
\end{equation}
and%
\begin{equation}
f_{k}(z)=-%
{\displaystyle\sum\limits_{j=1}^{k}}
\frac{1}{\left(  j+1\right)  !}\frac{d^{j}}{dz^{j}}f_{k-j}(z),\quad k\geq1,
\label{fk1}%
\end{equation}
where we have omitted (for the time being) any constant of integration. To
find the functions $f_{k}(z),$ we set
\begin{equation}
f_{k}(z)=a_{k}\frac{d^{k}}{dz^{k}}f_{0}(z) \label{fk2}%
\end{equation}
in (\ref{fk1}) and get $a_{0}=1$ and%
\[
a_{k}\frac{d^{k}}{dz^{k}}f_{0}(z)=-%
{\displaystyle\sum\limits_{j=1}^{k}}
\frac{1}{\left(  j+1\right)  !}\frac{d^{j}}{dz^{j}}\frac{d^{k-j}}{dz^{k-j}%
}f_{0}(z),\quad k\geq1,
\]
which gives%
\[
a_{0}=1,\quad a_{k}=-%
{\displaystyle\sum\limits_{j=1}^{k}}
\frac{1}{\left(  j+1\right)  !}a_{k-j},\quad k\geq1,
\]
or%
\begin{equation}
a_{0}=1,\quad%
{\displaystyle\sum\limits_{j=0}^{k}}
\frac{1}{\left(  k+1-j\right)  !}a_{j}=0,\quad k\geq1. \label{a1}%
\end{equation}
Multiplying both sides by $(k+1)!,$ we can write (\ref{a1}) as%
\[
a_{0}=1,\quad%
{\displaystyle\sum\limits_{j=0}^{k}}
\frac{(k+1)!}{\left(  k+1-j\right)  !j!}j!a_{j}=0,\quad k\geq1,
\]
or%
\begin{equation}
a_{0}=1,\quad%
{\displaystyle\sum\limits_{j=0}^{k}}
\binom{k+1}{j}j!a_{j}=0,\quad k\geq1. \label{a2}%
\end{equation}
Comparing (\ref{bernoulli}) and (\ref{a2}) we conclude that%
\begin{equation}
a_{k}=\frac{B_{k}}{k!},\quad k\geq0. \label{ak}%
\end{equation}
Thus, from (\ref{f0}), (\ref{fk2}) and (\ref{ak}) we have%
\[
f_{k}(z)=\frac{B_{k}}{k!}\frac{d^{k}}{dz^{k}}\left[  z\ln(z)-z\right]  ,\quad
k\geq0,
\]
from which we obtain%
\begin{equation}
f_{1}(z)=-\frac{1}{2}\ln(z) \label{f1}%
\end{equation}
and%
\[
f_{k}(z)=\frac{B_{k}}{k!}\left(  -1\right)  ^{k}\frac{\left(  k-1\right)
!}{z^{k-1}}=\frac{\left(  -1\right)  ^{k}B_{k}}{k\left(  k+1\right)  z^{k-1}%
},\quad k\geq2.
\]
Since $B_{2k+1}=0$ for all $k\geq1,$ we need to consider even values of $k$
only,%
\begin{equation}
f_{2k}(z)=\frac{B_{2k}}{2k\left(  2k+1\right)  z^{2k-1}},\quad k\geq1.
\label{f2}%
\end{equation}

So far, we haven't included any constant of integration in our calculations.
We could add a constant to one of the functions $f_{k}(z),$ let's say to
$f_{1}(z),$ and proceed as in Section \ref{section} to find it. Doing this, we
would obtain from (\ref{f0}), (\ref{f1}) and (\ref{f2}) that%
\[%
{\displaystyle\sum\limits_{k=0}^{N}}
f_{k}(z)=\mathcal{P}(z)+R_{N}(z)
\]
where $\mathcal{P}(z),$ $R_{N}(z)$ were defined in (\ref{P}) and (\ref{R}) respectively.

Another possibility, would be to assume no previous knowledge of $\Gamma(z),$
except for the difference equation $\Gamma(z+1)=z\Gamma(z)$ and the value at
$1,$ $\Gamma(1)=1.$ In doing so, (\ref{ser1}) would imply the initial
conditions $f_{k}(1)=0,$ for all $k\geq0.$ Hence, we would have%
\begin{align*}
f_{0}(z) &  =z\ln(z)-z+1,\quad f_{1}(z)=-\frac{1}{2}\ln(z)\\
f_{2k}(z) &  =\frac{B_{2k}}{2k\left(  2k+1\right)  z^{2k-1}}-\frac{B_{2k}%
}{2k\left(  2k+1\right)  },\quad k\geq1
\end{align*}
and therefore%
\begin{equation}
\ln\left[  \Gamma(z)\right]  \sim z\ln(z)-z-\frac{1}{2}\ln(z)+C_{N}%
+R_{N}(z),\label{ser4}%
\end{equation}
with%
\begin{equation}
C_{N}=1-%
{\displaystyle\sum\limits_{k=1}^{N}}
\frac{B_{2k}}{2k\left(  2k+1\right)  }.\label{Cn}%
\end{equation}
Computing the first few $C_{N}$, we would get%
\begin{align*}
C_{1} &  \simeq.91667,\ C_{2}\simeq.91944,\ C_{3}\simeq.91865,\ C_{4}%
\simeq.91925,\ C_{5}\simeq.91840,\\
C_{6} &  \simeq.92032,\ C_{7}\simeq.91391,\ C_{8}\simeq.94346,\ C_{9}%
\simeq.76382,\ C_{10}\simeq2.1562,
\end{align*}
and increasingly greater numbers (in absolute value). We would conclude that,
before the sum starts diverging, the $C_{N}$ seem to be approaching a value
close to $.918.$ Given our previous discussion of Stirling's formula, it is
not surprising to find that $\frac{1}{2}\ln\left(  2\pi\right)  \simeq.91894.$
Of course, geniuses like Euler or Gauss would reach the conclusion that the
optimal constant equals $\frac{1}{2}\ln\left(  2\pi\right)  $ without knowing
anything about Stirling's work!

\section{Conclusion}

We have presented the history and previous approaches to the proof of
Stirling's series (\ref{Gamma2}). We have used a different procedure, based on
the asymptotic analysis of the difference equation (\ref{difference}). The
method reproduces (\ref{Gamma2}) very easily and can be extended to use in
more complicated difference equations.

Bender and Orszag observed in \cite[Page 227]{MR538168} that

\begin{quotation}
without further information the constant $\frac{1}{2}\ln\left(  2\pi\right)  $
cannot be determined. The difference equation that we have solved is linear
and homogeneous, so any arbitrary multiple of a solution is still a solution.
\end{quotation}

While agreeing with them completely, we have shown that by imposing the
additional condition $\Gamma(1)=1$ one can find an approximation to the value
of $\frac{1}{2}\ln\left(  2\pi\right)  ,$ without any other assumptions. Thus,
local behavior at $z=1$ and asymptotic behavior as $z\rightarrow\infty$ can be
combined fruitfully.

We sincerely hope that more and more professors and their students will
discover the beauty contained in the (very!) Special Functions, among which
$\Gamma(z)$ is, without doubt, a prima donna.

\strut

\end{document}